\documentclass[leqno,12pt]{article}
\usepackage{}

\usepackage{amssymb}
\usepackage{euscript}
\usepackage[dvips]{graphicx}
\usepackage{flafter}
\usepackage{pstricks}
\usepackage[latin1]{inputenc}
\usepackage{amsmath}
\usepackage{amsfonts}
\usepackage{enumitem}

\usepackage{verbatim}

\usepackage{mathrsfs} 

\evensidemargin\oddsidemargin

\begin{document}

\baselineskip=18pt
\setcounter{page}{1}

\renewcommand{\theequation}{\thesection.\arabic{equation}}
\newtheorem{theorem}{Theorem}[section]
\newtheorem{lemma}[theorem]{Lemma}
\newtheorem{proposition}[theorem]{Proposition}
\newtheorem{corollary}[theorem]{Corollary}
\newtheorem{remark}[theorem]{Remark}
\newtheorem{fact}[theorem]{Fact}
\newtheorem{problem}[theorem]{Problem}
\newtheorem{example}[theorem]{Example}
\newtheorem{question}[theorem]{Question}
\newtheorem{conjecture}[theorem]{Conjecture}

\newcommand{\eqnsection}{
\renewcommand{\theequation}{\thesection.\arabic{equation}}
    \makeatletter
    \csname  @addtoreset\endcsname{equation}{section}
    \makeatother}
\eqnsection

\def\r{{\mathbb R}}
\def\e{{\mathbb E}}
\def\p{{\mathbb P}}
\def\P{{\bf P}}
\def\E{{\bf E}}
\def\Q{{\bf Q}}
\def\z{{\mathbb Z}}
\def\N{{\mathbb N}}
\def\T{{\mathbb T}}
\def\G{{\mathbb G}}
\def\L{{\mathbb L}}

\def\deg{\chi}

\def\ee{\mathrm{e}}
\def\d{\, \mathrm{d}}
\def\S{\mathscr{S}}



\vglue50pt

\centerline{\large\bf  Scaling limit of the path leading to the leftmost particle}
\centerline{\large\bf  in a branching random walk}

\bigskip
\bigskip

\centerline{Xinxin CHEN}

\medskip

\centerline{\it Universit\'e Paris VI}

\bigskip
\bigskip
\bigskip

{\leftskip=2truecm \rightskip=2truecm \baselineskip=15pt \small

\noindent{\slshape\bfseries Summary.} We consider a discrete-time branching random walk defined on the real line, which is assumed to be supercritical and in the boundary case. It is known that its leftmost position of the $n$-th generation behaves asymptotically like $\frac{3}{2}\ln n$, provided the non-extinction of the system. The main goal of this paper, is to prove that the path from the root to the leftmost particle, after a suitable normalizatoin, converges weakly to a Brownian excursion in $D([0,1],\r)$.

\bigskip

\bigskip

} 

\bigskip
\bigskip

\section{Introduction}
   \label{s:intro}

$\phantom{aob}$We consider a branching random walk, which is constructed according to a point process $\mathcal{L}$ on the line. Precisely speaking, the system is started with one initial particle at the origin. This particle is called the root, denoted by $\varnothing$. At time $1$, the root dies and gives birth to some new particles, which form the first generation. Their positions constitute a point process distributed as $\mathcal{L}$. At time 2, each of these particles dies and gives birth to new particles whose positions -- relative to that of their parent -- constitute a new independent copy of $\mathcal{L}$. The system grows according to the same mechanism.

We denote by $\mathbb{T}$ the genealogical tree of the system, which is clearly a Galton-Watson tree rooted at $\varnothing$. If a vertex $u\in\mathbb{T}$ is in the $n$-th generation, we write $|u|=n$ and denote its position by $V(u)$. Then $\{V(u),\; |u|=1\}$ follows the same law as $\mathcal{L}$. The family of positions $(V(u);\; u\in\mathbb{T})$ is viewed as our branching random walk.

Throughout the paper, the branching random walk is assumed to be in the boundary case (Biggins and Kyprianou~\cite{biggins-kyprianou05}):
\begin{equation}
    \E\Big[\sum_{|u|=1}1\Big]>1,
    \qquad
    \E \Big[ \sum_{|x|=1} \ee^{-V(x)} \Big]=1,
    \qquad
    \E \Big[ \sum_{|x|=1} V(x) \ee^{-V(x)} \Big] =0.
    \label{cond-boundary}
\end{equation}
For any $y\in\r$, let $y_+:=\max\{y,\, 0\}$ and $\log_+y:=\log(\max\{y,\, 1\})$. We also assume the following integrability conditions:
\begin{eqnarray}
    \label{cond-2}
    \E\Big[\sum_{|u|=1} V(u)^2\ee^{-V(u)}\Big]
 &<&\infty,
    \\
    \label{cond-square}
    \E [X (\log_+ X)^2 ]
 &<&\infty,
    \qquad
    \E[\widetilde{X} \log_+ \widetilde{X}] <\infty,
    \label{cond-x}
\end{eqnarray}

\noindent where
$$
    X:=\sum_{|u|=1}\ee^{-V(u)},
    \qquad
    \widetilde{X}
    :=
    \sum_{|u|=1} V(u)_+ \ee^{-V(u)}.
$$
We define $I_n$ to be the leftmost position in the $n$-th generation, i.e.
\begin{equation}
I_n:=\inf\{V(u),\; |u|=n\},
\label{def-leftmostposition}
\end{equation}
with $\inf\emptyset:=\infty$. If $I_n<\infty$, we choose a vertex uniformly in the set $\{u: |u|=n,\; V(u)=I_n\}$ of leftmost particles at time $n$ and denote it by $m^{(n)}$. We let $[\![ \varnothing, \, m^{(n)}]\!]=\{\varnothing=:m_0^{(n)},\, m_1^{(n)},\dots,\, m_n^{(n)}:=m^{(n)}\}$ be the shortest path in $\mathbb{T}$ relating the root $\varnothing$ to $m^{(n)}$, and introduce the path from the root to $m^{(n)}$ as follows
\begin{equation*}
(I_n(k);\; 0\leq k\leq n):= (V(m^{(n)}_{k});\; 0\leq k\leq n).
\end{equation*}
In particular, $I_n(0)=0$ and $I_n(n)=I_n$. Let $\sigma$ be the positive real number such that $\sigma^2=\E\Big[\sum_{|u|=1} V(u)^2\ee^{-V(u)}\Big]$. Our main result is as follows.
\begin{theorem}\label{mainthm}
The rescaled path $(\frac{I_n(\lfloor sn\rfloor)}{\sigma\sqrt{n}};\; 0\leq s\leq 1)$ converges in law in $D([0,1],\r)$, to a normalized Brownian excursion $(e_s;\; 0\leq s\leq 1)$.
\end{theorem}
\begin{remark}
It has been proved in \cite{addario-berry-reed}, \cite{hu-shi} and \cite{elie} that $I_n$ is around $\frac{3}{2}\ln n$. In \cite{elie-berestycki-brunet-shi}, the authors proved that, for the model of branching Brownian motion, the time reversed path followed by the leftmost particle converges in law to a certain stochastic process.
\end{remark}

Let us say a few words about the proof of Theorem \ref{mainthm}. We first consider the path leading to $m^{(n)}$, by conditioning that its ending point $I_n$ is located atypically below $\frac{3}{2}\ln n-z$ with large $z$. Then we apply the well-known spinal decomposition to show that this path, conditioned to $\{I_n\leq \frac{3}{2}\ln n-z\}$, behaves like a simple random walk staying positive but tied down at the end. Such a random walk, being rescaled, converges in law to the Brownian excursion (see \cite{durrett}). We then prove our main result by removing the condition of $I_n$.  The main strategy is borrowed from \cite{elie}, but with appropriate refinements.

The rest of the paper is organized as follows. In Section 2, we recall the spinal decomposition by a change of measures, which implies the useful many-to-one lemma. We prove a conditioned version of Theorem \ref{mainthm} in Section 3. In Section 4, we remove the conditioning and prove the theorem.

Throughout the paper, we use $a_n \sim b_n$ ($n\to \infty$) to denote $\lim_{n\to \infty} \, {a_n\over b_n} =1$; and let $(c_i)_{i\geq 0}$ denote finite and positive constants. We write $\E[f;\,A]$ for $\E[f\textbf{1}_A]$. Moreover, $\sum_\varnothing := 0$ and $\prod_\varnothing :=1$.

\section{Lyons' change of measures and spinal decomposition}
\label{s:new-pair}

$\phantom{aob}$ For any $a\in\r$, let $\P_a$ be the probability measure such that $\P_a((V(u),\; u\in\mathbb{T})\in\cdot)=\P((V(u)+a,\; u\in\mathbb{T})\in\cdot)$. The corresponding expectation is denoted by $\E_a$. Let $(\mathcal{F}_n,\; n\geq 0)$ be the natural filtration generated by the branching random walk and let $\mathcal{F}_\infty:=\vee_{n\geq 0}\mathcal{F}_n$. We introduce the following random variables:
\begin{equation}
W_n:=\sum_{|u|=n}e^{-V(u)},
\qquad n\geq 0.
\label{addmart}
\end{equation}
It follows immediately from (\ref{cond-boundary}) that $(W_n,\; n\geq 0)$ is a non-negative martingale with respect to $(\mathcal{F}_n)$. It is usually referred as the additive martingale. We define a probability measure $\Q_a$ on $\mathcal{F}_\infty$ such that for any $n\geq 0$,
\begin{equation}\label{RNdensity}
\frac{d\Q_a}{d\P_a}\bigg\vert_{\mathcal{F}_n}:= e^a W_n.
\end{equation}
For convenience, we write $\Q$ for $\Q_0$.

 Let us give the description of the branching random walk under $\Q_a$ in an intuitive way, which is known as the spinal decomposition. We introduce another point process $\widehat{\mathcal{L}}$ with Radon-Nykodin derivative $\sum_{x\in\mathcal{L}}e^{-x}$ with respect to the law of $\mathcal{L}$. Under $\Q_a$, the branching random walk evolves as follows. Initially, there is one particle $w_0$ located at $V(w_0)=a$. At each step $n$, particles at generation $n$ die and give birth to new particles independently according to the law of $\mathcal{L}$, except for the particle $w_n$ which generates its children according to the law of $\widehat{\mathcal{L}}$. The particle $w_{n+1}$ is chosen proportionally to $e^{-V(u)}$ among the children $u$ of $w_n$. We still call $\mathbb{T}$ the genealogical tree of the process, so that $(w_n)_{n\geq0}$ is a ray in $\mathbb{T}$, which is called the spine. This change of probabilities was presented in various forms; see, for example \cite{lyons}, \cite{hu-shi} and \cite{chauvin-rouault}.

It is convenient to use the following notation. For any $u\in\mathbb{T}\setminus\{\varnothing\}$, let $\overleftarrow{u}$ be the parent of $u$, and
\begin{equation*}
\Delta V(u):=V(u)-V(\overleftarrow{u}).
\end{equation*}
Let $\Omega(u)$ be the set of brothers of $u$, i.e. $\Omega(u):=\{v\in\mathbb{T}: \overleftarrow{v}=\overleftarrow{u},\; v\neq u\}$. Let $\delta$ denote the Dirac measure. Then under $\Q_a$, $\sum_{|u|=1}\delta_{\Delta V(u)}$ follows the law of $\widehat{\mathcal{L}}$. Further, We recall the following proposition, from \cite{hu-shi} and \cite{lyons}.
\begin{proposition}\label{changeofmeasure}
\begin{itemize}[fullwidth]
  \item [(1)]For any $|u|=n$, we have
  \begin{equation}\Q_a[w_n=u\vert\mathcal{F}_n]=\frac{e^{-V(u)}}{W_n}.\end{equation}
  \item [(2)]Under $\Q_a$, the random variables $\Big(\sum_{v\in\Omega(w_n)}\delta_{\Delta V(v)},\; \Delta V(w_n)\Big),\; n\geq 1$ are i.i.d..
\end{itemize}
\end{proposition}

As a consequence of this proposition, we get the many-to-one lemma as follows:
\begin{lemma}
There exists a centered random walk $(S_n;\; n\geq 0)$ with $\P_a(S_0=a)=1$ such that for any $n\geq 1$ and any measurable function $g:\mathbb{R}^n\rightarrow[0,\infty)$, we have
\begin{equation}\label{many-to-one}
\E_a\bigg[\sum_{|u|=n}g(V(u_1),\dots,V(u_n))\bigg]=\E_a[e^{S_n-a}g(S_1,\dots, S_n)],
\end{equation}
where we denote by $[\![ \varnothing, \, u]\!]=\{\varnothing=:u_0, \, u_1\dots, u_{|u|}:=u\}$ the ancestral line of $u$ in $\mathbb{T}$.
\end{lemma}

Note that by (\ref{cond-square}), $S_1$ has the finite variance $\sigma^2=\E[S_1^2]=\E[\sum_{|u|=1}V(u)^2e^{-V(u)}]$.

\subsection{Convergence in law for the one-dimensional random walk}
Let us introduce some results about the centered random walk $(S_n)$ with finite variance, which will be used later. For any $0\leq m\leq n$, we define $\underline{S}_{[m,n]}:=\min_{m\leq j\leq n}S_j$, and $\underline{S}_n=\underline{S}_{[0,n]}$. We denote by $R(x)$ the renewal function of $(S_n)$, which is defined as follows:
\begin{equation}
R(x)=\textbf{1}_{\{x=0\}}+\textbf{1}_{\{x>0\}}\sum_{k\geq0}\P(-x\leq S_k<\underline{S}_{n-1}).
\end{equation}
For the random walk $(-S_n)$, we define $\underline{S}^-_{[m,n]}$, $\underline{S}^-_n$ and $R_-(x)$ similarly. It is known (see \cite{feller} p. 360) that there exists $c_0>0$ such that
\begin{equation}
\lim_{x\rightarrow\infty}\frac{R(x)}{x}=c_0.
\end{equation}
Moreover, it is shown in \cite{kozlov} that there exist $C_+,\; C_->0$ such that for any $a\geq 0$,
\begin{align}
\P_a\Big(\underline{S}_n\geq 0\Big)&\sim \frac{C_+}{\sqrt{n}}R(a);\label{probstaypositive}\\
\P_{a}\Big(\underline{S}^-_n\geq0\Big)&\sim \frac{C_-}{\sqrt{n}}R_-(a)\label{probstaynegative}.
\end{align}
We also state the following inequalities (see Lemmas 2.2 and 2.4 in \cite{elie-shi1}, respectively).
\begin{fact}\label{basic}
\begin{description}[fullwidth]
\item[ (i)]  There exists a constant $c_1>0$ such that for any $b\geq a\geq 0$, $x\geq0$ and $n\geq 1$,
\begin{equation}\label{element}
\P\big(\underline{S}_n\geq -x;\; S_n\in[a-x, b-x]\big)\leq c_1(1+x)(1+b-a)(1+b)n^{-3/2}.
\end{equation}
\item[ (ii)] Let $0<\lambda<1$. There exists a constant $c_2>0$ such that for any $b\geq a\geq 0$, $x,\; y\geq0$ and $n\geq 1$,
  \begin{equation}\label{elementary}
  \P_x(S_n\in[y+a,y+b],\underline{S}_n\geq0, \underline{S}_{[\lambda n, n]}\geq y) \leq c_2(1+x)(1+b-a)(1+b)n^{-3/2}.
  \end{equation}
\end{description}
\end{fact}

{\color{black}Before we give the next lemma, we recall the definition of lattice distribution (see \cite{feller}, p. 138). The distribution of a random variable $X_1$ is lattice, if it is concentrated on a set of points $\alpha+\beta \z$, with $\alpha$ arbitrary. The largest $\beta$ satisfying this property is called the span of $X_1$. Otherwise, the distribution of $X_1$ is called non-lattice. }


\begin{lemma}\label{keyconv}
Let $(r_n)_{n\geq 0}$ be a sequence of real numbers such that $\lim_{n\rightarrow \infty}\frac{r_n}{\sqrt{n}}=0$. Let $f:\mathbb{R}_+\rightarrow \r$ be a Riemann integrable function. We suppose that there exists a non-increasing function $\overline{f}:\r_+\rightarrow\r$ such that $|f(x)|\leq \overline{f}(x)$ for any $x\geq 0$ and $\int_{x\geq 0}x\overline{f}(x)dx<\infty$. For $0<\Delta<1$, let $F:D([0,\Delta],\; \r)\rightarrow [0,1]$ be continuous. Let $a\geq0$.
\begin{description}[fullwidth]
  \item[(I) Non-lattice case.] If the distribution of $(S_1-S_0)$ is non-lattice, then there exists a constant $C_1>0$ such that
  \begin{multline}\label{convwithkilling}
  \lim_{n\rightarrow\infty}n^{3/2}\E\Big[F\Big(\frac{S_{\lfloor sn\rfloor}}{\sigma\sqrt{n}}; 0\leq s\leq \Delta\Big)f(S_n-y);\;\underline{S}_n\geq -a,\, \underline{S}_{[\Delta n, n]}\geq y\Big]\\
  =C_1R(a)\int_{x\geq 0}f(x)R_-(x)dx\E[F(e_s;0\leq s\leq \Delta)],
  \end{multline}
  uniformly in $y\in[0,r_n]$.
  \item[(II) Lattice case.] If the distribution of $(S_1-S_0)$ is supported in $(\alpha+\beta\z)$ with span $\beta$, then for any $d\in \r$,
      \begin{multline}\label{convwithkilling-lattice}
      \lim_{n\rightarrow\infty}n^{3/2}\E\Big[F\Big(\frac{S_{\lfloor sn\rfloor}}{\sigma\sqrt{n}}; 0\leq s\leq \Delta\Big)f(S_n-y+d);\;\underline{S}_n\geq -a,\,\underline{S}_{[\Delta n, n]}\geq y-d\Big]\\
       =C_1R(a)\beta\sum_{j\geq \lceil-\frac{d}{\beta}\rceil}f(\beta j+d)R_-(\beta j+d)\E[F(e_s;0\leq s\leq \Delta)].
      \end{multline}
  uniformly in $y\in[0,r_n]\cap\{\alpha n+\beta\z\}$.
\end{description}
\end{lemma}

\textit{Proof of Lemma \ref{keyconv}.} The lemma is a refinement of Lemma 2.3 in \cite{elie}, which proved the convergence in the non-lattice case when $a=0$ and $F\equiv1$. We consider the non-lattice case first. We denote the expectation on the left-hand side of (\ref{convwithkilling}) by $\chi(F,f)$. Observe that for any $K\in\mathbb{N}_+$,
\begin{equation*}
\chi(F,f)=\chi\big(F, f(x)1_{(0\leq x\leq K)}\big)+\chi\big(F,f(x)1_{(x>K)}\big).
\end{equation*}
Since $0\leq F\leq 1$, we have $\chi\big(F,f(x)1_{(x>K)}\big) \leq  \chi\big(1,f(x)1_{(x>K)}\big)$, which is bounded by
\begin{equation*}
\sum_{j\geq K}\E_a\Big[f(S_n-y-a);\;\underline{S}_n\geq 0,\, \underline{S}_{[\Delta n, n]}\geq y+a,\, S_n\in[y+a+j, y+a+j+1]\Big].
\end{equation*}
Recall that $|f(x)|\leq \overline{f}(x)$ with $\overline{f}$ non-increasing. We get that
\begin{equation*}
\chi\big(1,f(x)1_{(x>K)}\big) \leq \sum_{j\geq K}\overline{f}(j)\P_a\Big[\underline{S}_n\geq 0,\, \underline{S}_{[\Delta n, n]}\geq y+a,\, S_n\in[y+a+j, y+a+j+1]\Big].
\end{equation*}
It then follows from (\ref{elementary}) that
\begin{equation}
\chi\big(1,f(x)1_{(x>K)}\big) \leq 2c_2(1+a)\Big(\sum_{j\geq K}\overline{f}(j)(2+j)\Big)n^{-3/2}.
\end{equation}
Since $\int_{0}^\infty x\overline{f}(x)dx<\infty$, the sum $\sum_{j\geq K}\overline{f}(j)(2+j)$ decreases to zero as $K\uparrow \infty$. We thus only need to estimate $\chi\big(F, f(x)1_{(0\leq x\leq K)}\big)$. Note that $f$ is Riemann integrable. It suffices to consider $\chi\Big(F, 1_{(0\leq x\leq K)}\Big)$ with $K$ a positive constant.

Applying the Markov property at time $\lfloor\Delta n\rfloor$ shows that
\begin{eqnarray}\label{markovprop}
\chi\Big(F, 1_{(0\leq x\leq K)}\Big)&=&\E_a\Big[F\Big(\frac{S_{\lfloor sn\rfloor}-a}{\sigma\sqrt{n}}; 0\leq s\leq \Delta\Big);\; S_n\leq y+a+K,\underline{S}_n\geq 0, \underline{S}_{[\Delta n, n]}\geq y+a\Big]\nonumber\\
&=&\E_a\Big[F\Big(\frac{S_{\lfloor sn\rfloor}-a}{\sigma\sqrt{n}}; 0\leq s\leq \Delta\Big)\Psi_K(S_{\lfloor\Delta n\rfloor});\; \underline{S}_{ \lfloor\Delta n\rfloor}\geq 0\Big],
\end{eqnarray}
where $\Psi_K(x):=\P_x\Big[S_{n-\lfloor\Delta n\rfloor}\leq y+a+K,\; \underline{S}_{n-\lfloor\Delta n\rfloor}\geq y+a\Big]$. By reversing time, we obtain that $\Psi_K(x)=\P\Big[\underline{S}^-_m\geq (-S_m)+(y+a-x)\geq-K\Big]$ with $m:=n-\lfloor\Delta n\rfloor$.

We define $\tau_n$ as the first time when the random walk $(-S)$ hits the minimal level during $[0,n]$, namely, $\tau_n:=\inf\{k\in[0,n]: -S_k=\underline{S}^-_n\}$. Define also $\varkappa(z,\zeta; n):=\P(-S_n\in[z,z+\zeta],\;\underline{S}^-_n\geq 0)$ for any $z$, $\zeta\geq 0$. Then,
\begin{equation}\label{hittingmin}
\begin{split}
\Psi_K(x)&= \sum_{k=0}^m \P\Big[\tau_m=k;\; \underline{S}^-_m\geq (-S_m)+(y+a-x)\geq-K\Big]\\
         &= \sum_{k=0}^m \P\Big[-S_k=\underline{S}^-_k\geq -K;\;\varkappa(x-y-a, \underline{S}^-_k+K; m-k)\Big],
\end{split}
\end{equation}
where the last equality follows from the Markov property.

Let $\psi(x):=xe^{-x^2/2}\textbf{1}_{(x\geq0)}$. Combining Theorem 1 of \cite{caravenna} with (\ref{probstaypositive}) yields that
\begin{equation}\label{locallim}
\varkappa(z,\zeta; n)=\P_0\Big[-S_n\in[z,z+\zeta];\; \underline{S}_n\geq 0\Big]=\frac{C_-\zeta}{\sigma n}\psi\Big(\frac{z}{\sigma\sqrt{n}}\Big)+o(n^{-1}),
\end{equation}
uniformly in $z\in\r_+$ and $\zeta$ in compact sets of $\r_+$. Note that $\psi$ is bounded on $\r_+$. Therefore, there exists a constant $c_3>0$ such that for any $\zeta\in[0,K]$, $z\geq 0$ and $n\geq 0$,
\begin{equation}\label{globalupp}
\varkappa(z,\zeta; n)\leq c_3 \frac{(1+K)}{n+1}.
\end{equation}

Let $k_n:=\lfloor \sqrt{n} \rfloor$. We divide the sum on the right-hand side of (\ref{hittingmin}) into two parts:
\begin{equation}\label{division}
\Psi_K(x)=\sum_{k=0}^{k_n}+\sum_{k=k_n+1}^m \P\big[-S_k=\underline{S}^-_k\geq -K;\;\varkappa(x-y-a, \underline{S}^-_k+K; m-k)\big].
\end{equation}
By (\ref{locallim}), under the assumption that $y=o(\sqrt{n})$, the first part becomes that
\begin{eqnarray}\label{partone}
&&\frac{C_-}{\sigma m}\psi\Big(\frac{x-a}{\sigma\sqrt{m}}\Big)\sum_{k=0}^{k_n}\E\big[\underline{S}^-_k+K; -S_k=\underline{S}^-_k\geq -K\big]
+o(n^{-1})\sum_{k=0}^{k_n}\P\big[-S_k=\underline{S}^-_k\geq -K\big]\\
&=&\frac{C_-}{\sigma m}\psi\Big(\frac{x-a}{\sigma\sqrt{m}}\Big)\int_0^K R_-(u)du+o(n^{-1}),\nonumber
\end{eqnarray}
where the last equation comes from the fact that $\sum_{k\geq0}\E\big[\underline{S}^-_k+K; -S_k=\underline{S}^-_k\geq -K\big]=\int_{0}^K R_-(u)du$.
On the other hand, using (\ref{globalupp}) for $\varkappa(x-y-a, \underline{S}^-_k+K; m-k)$ and then applying (i) of Fact \ref{basic} imply that for $n$ large enough, the second part of (\ref{division}) is bounded by
\begin{equation}\label{parttwo}
\begin{split}
& \sum_{k=k_n+1}^m c_3\frac{1+K}{m+1-k}\P\big(\underline{S}^-_k\geq -K,\; -S_k\in[-K,0]\big)\\
&\leq c_4\sum_{k=k_n+1}^m \frac{(1+K)^3}{(m+1-k)k^{3/2}}=o(n^{-1}).
\end{split}
\end{equation}
By (\ref{partone}) and (\ref{parttwo}), we obtain that as $n$ goes to infinity,
\begin{equation}
\Psi_K(x)=o(n^{-1})+\frac{C_-}{\sigma(n-\lfloor\Delta n\rfloor)}\psi\Big(\frac{x-a}{\sigma\sqrt{n-\lfloor\Delta n\rfloor}}\Big)\int_{0}^K R_-(u)du,
\end{equation}
uniformly in $x\geq 0$ and $y\in[0,r_n]$. Plugging it into (\ref{markovprop}) and then combining with (\ref{probstaypositive}) yield that
\begin{eqnarray}
\chi(F, 1_{(0\leq x\leq K)})&=&o(n^{-3/2})+\frac{C_-}{\sigma(1-\Delta)n}\int_{0}^K R_-(u)du\nonumber\\
&&\times\frac{C_+R(a)}{\sqrt{\Delta n}}\E_a\Big[F\Big(\frac{S_{\lfloor sn\rfloor}-a}{\sigma\sqrt{n}}; 0\leq s\leq \Delta\Big)\psi\Big(\frac{S_{\Delta n}-a}{\sigma\sqrt{(1-\Delta)n}}\Big)\Big\vert\underline{S}_{\Delta n}\geq0\Big].\nonumber
\end{eqnarray}

Theorem 1.1 of \cite{caravenna-loic} says that under the conditioned probability $\P_a\Big(\cdot \Big\vert \underline{S}_{ \Delta n}\geq0\Big)$, $(\frac{S_{\lfloor r\Delta n\rfloor}}{\sigma\sqrt{\Delta n}};0\leq r\leq 1)$ converges in law to a Brownian meander, denoted by $(\mathcal{M}_{r};0\leq r\leq 1)$. Therefore,
\begin{multline*}
\chi(F,1_{(0\leq x\leq K)})\sim \frac{C_-C_+R(a)}{\sigma n^{3/2}(1-\Delta)\sqrt{\Delta}}\int_{0}^K R_-(u)du\E\Big[F\Big(\sqrt{\Delta}\mathcal{M}_{s/\Delta};0\leq s\leq\Delta\Big)\psi\Big(\frac{\sqrt{\Delta}\mathcal{M}_1}{\sqrt{1-\Delta}}\Big)\Big].
\end{multline*}
It remains to check that
\begin{equation}\label{meandertoex}
\frac{1}{(1-\Delta)\sqrt{\Delta}}\E\Big[F\Big(\sqrt{\Delta}\mathcal{M}_{s/\Delta};0\leq s\leq\Delta\Big)\psi\Big(\frac{\sqrt{\Delta}\mathcal{M}_1}{\sqrt{1-\Delta}}\Big)\Big]=\sqrt{\frac{\pi}{2}}\E\Big[F(e_s;0\leq s\leq \Delta)\Big].
\end{equation}
Let $(R_s; 0\leq s\leq1)$ be a standard three-dimensional Bessel process. Then, as is shown in \cite{imhof},
\begin{eqnarray*}
&&\frac{1}{(1-\Delta)\sqrt{\Delta}}\E\Big[F\Big(\sqrt{\Delta}\mathcal{M}_{s/\Delta};0\leq s\leq\Delta\Big)\psi\Big(\frac{\sqrt{\Delta}\mathcal{M}_1}{\sqrt{1-\Delta}}\Big)\Big]\\
&=&\sqrt{\frac{\pi}{2}}\frac{1}{(1-\Delta)\sqrt{\Delta}}\E\Big[\frac{1}{R_1}F\Big(\sqrt{\Delta}R_{s/\Delta};0\leq s\leq\Delta\Big)\psi\Big(\frac{\sqrt{\Delta}R_1}{\sqrt{1-\Delta}}\Big)\Big],\\
&=&\sqrt{\frac{\pi}{2}}\E\Big[\frac{1}{(1-\Delta)^{3/2}}e^{-\frac{R_\Delta^2}{2(1-\Delta)}}F\Big(R_{s};0\leq s\leq\Delta\Big)\Big],
\end{eqnarray*}
where the last equation follows from the scaling property of Bessel process. Let $(r_s; 0\leq s\leq 1)$ be a standard three-dimensional Bessel bridge. Note that for any $\Delta<1$, $(r_s;0\leq s\leq \Delta)$ is equivalent to $(R_s; 0\leq s\leq \Delta)$, with density $(1-\Delta)^{-3/2}\exp(-\frac{R_\Delta^2}{2(1-\Delta)})$ (see p. 468 (3.11) of \cite{yor}). Thus,
\begin{equation*}
\frac{1}{(1-\Delta)\sqrt{\Delta}}\E\Big[F\Big(\sqrt{\Delta}\mathcal{M}_{s/\Delta};0\leq s\leq\Delta\Big)\psi\Big(\frac{\sqrt{\Delta}\mathcal{M}_1}{\sqrt{1-\Delta}}\Big)\Big]=\sqrt{\frac{\pi}{2}}\E\Big[F(r_s; 0\leq s\leq\Delta)\Big].
\end{equation*}
Since a normalized Brownian excursion is exactly a standard three-dimensional Bessel bridge, this yields (\ref{meandertoex}). Therefore, (\ref{convwithkilling}) is proved by taking $C_1=\sqrt{\frac{\pi}{2}}\frac{C_-C_+}{\sigma}$.

The proof of the lemma in the lattice case is along the same lines, except that we use Theorem 2 (instead of Theorem 1) of \cite{caravenna}.    $\square$

 \section{Conditioning on the event $\{I_n\leq \frac{3}{2}\ln n-z\}$}

 $\phantom{aob}$On the event $\{I_n\leq \frac{3}{2}\ln n-z\}$, we analyze the sample path leading to a particle located at the leftmost position at the $n$th generation. For $z\geq 0$ and $n\geq1$, let $a_n(z):=\frac{3}{2}\ln n-z$ if the distribution of $\mathcal{L}$ is non-lattice and let $a_n(z):=\alpha n+\beta\lfloor\frac{\frac{3}{2}\ln n-\alpha n}{\beta}\rfloor-z$ if the distribution of $\mathcal{L}$ is supported by $\alpha+\beta\z$.
 This section is devoted to the proof of the following proposition.
\begin{proposition}\label{conditionedconv}
For any $\Delta\in(0,1]$ and any continuous functional $F: D([0,\Delta],\; \r)\rightarrow[0,1]$,
\begin{equation}
\lim_{z\rightarrow\infty}\limsup_{n\rightarrow\infty}\bigg\vert\E\Big[F\Big(\frac{I_n(\lfloor sn \rfloor)}{\sigma\sqrt{n}};0\leq s\leq \Delta\Big)\Big\vert I_n\leq a_n(z) \Big]-\E\Big[F(e_s;0\leq s\leq \Delta)\Big]\bigg\vert=0.
\end{equation}
\end{proposition}

We begin with some preliminary results.

For any $0<\Delta<1$ and $L$, $K\geq 0$, we denote by $J_{z,K,L}^\Delta(n)$ the following collection of particles:
\begin{equation}
\Big\{u\in\mathbb{T}: |u|=n,\; V(u)\leq a_n(z),\; \min_{0\leq k\leq n}V(u_k)\geq-z+K,\;\min_{\Delta n\leq k\leq n}V(u_k)\geq a_n(z+L)\Big\}.
\end{equation}

\begin{lemma}\label{typicalmin}
For any $\varepsilon >0$, there exists $L_\varepsilon>0$ such that for any $L\geq L_\varepsilon$, $n\geq 1$ and $z\geq K\geq 0$,
\begin{equation}\label{pathbyleftmost}
\P\Big(m^{(n)}\not\in J^\Delta_{z,K,L}(n),\; I_n\leq a_n(z)\Big)\leq \Big(e^K+\varepsilon(1+z-K)\Big) e^{-z}.
\end{equation}
\end{lemma}

\noindent\textit{Proof.} It suffices to show that for any $\varepsilon\in(0,1)$, there exists $L_\varepsilon\geq 1$ such that for any $L\geq L_\varepsilon$, $n\geq 1$ and $z\geq K\geq 0$,
\begin{equation}
\P\Big(\exists |u|=n: V(u)\leq a_n(z),\; u\not\in J_{z,K,L}^\Delta(n)\Big)\leq \Big(e^K+\varepsilon(1+z-K)\Big) e^{-z}.
\end{equation}
We observe that
\begin{multline}\label{typicalcurve}
\P\Big(\exists |u|=n: V(u)\leq a_n(z),\; u\not\in J_{z,K,L}^\Delta(n)\Big)\leq \P\Big(\exists u\in\mathbb{T}: V(u)\leq -z+K \Big)\\
+\P\Big(\exists |u|=n: V(u)\leq a_n(z),\; \min_{0\leq k\leq n}V(u_k)\geq-z+K,\;\min_{\Delta n\leq k\leq n}V(u_k)\leq a_n(z+L)\Big).
\end{multline}
On the one hand, by (\ref{many-to-one}),
\begin{eqnarray}\label{firstbarrier}
\P\Big(\exists u\in\mathbb{T}: V(u)\leq -z+k \Big)&\leq& \sum_{n\geq 0}\E\bigg[\sum_{|u|=n}\textbf{1}_{\{V(u)\leq-z+K<\min_{k<n}V(u_k)\}}\bigg]\\
                                                  &=& \sum_{n\geq 0}\E[e^{S_n};\; S_n\leq -z+K<\underline{S}_{n-1}]\leq e^{-z+K}.\nonumber
\end{eqnarray}
On the other hand, denoting $A_n(z):=[a_n(z)-1, a_n(z)]$ for any $z\geq 0$,
\begin{eqnarray*}
&&\P\Big(\exists |u|=n: V(u)\leq a_n(z),\; \min_{0\leq k\leq n}V(u_k)\geq-z+K,\;\min_{\Delta n\leq k\leq n}V(u_k)\leq a_n(z+L)\Big)\\
&=&\P_{z-K}\Big(\exists |u|=n: V(u)\leq a_n(K),\;\min_{0\leq k\leq n}V(u_k)\geq0,\;\min_{\Delta n\leq k\leq n}V(u_k)\leq a_n(K+L)\Big)\\
&\leq&\sum_{\ell\geq L+K}\sum_{j=K}^{j=K+\ell}\P_{z-K}\Big(\exists |u|=n: V(u)\in A_n(j),\;\min_{0\leq k\leq n}V(u_k)\geq0,\;\min_{\Delta n\leq k\leq n}V(u_k)\in A_n(\ell)\Big).
\end{eqnarray*}
According to Lemma 3.3 in \cite{elie}, there exist constants $1>c_5>0$ and $c_6>0$ such that for any $n\geq 1$, $L\geq 0$ and $x$, $z\geq 0$,
\begin{eqnarray}\label{killingupp}
&&\P_x\Big(\exists u\in\mathbb{T}: |u|=n,\; V(u)\in A_n(z),\;\min_{0\leq k\leq n}V(u_k)\geq 0,\; \min_{\Delta n \leq k\leq n}V(u_k)\in A_n(z+L)\Big)\\
&\leq& c_6(1+x)e^{-c_5L}e^{-x-z}.\nonumber
\end{eqnarray}
Hence, combining (\ref{firstbarrier}) with (\ref{typicalcurve}) yields that
\begin{eqnarray*}
&&\P\Big(\exists |u|=n: V(u)\leq a_n(z),\; u\not\in J_{z,K,L}^\Delta(n)\Big)\\
&\leq& e^{-z+K}+ \sum_{\ell\geq L}\sum_{0\leq j\leq \ell}c_6(1+z-K)e^{-c_5(\ell-j)}e^{-z-j}\\
&\leq&\Big(e^K+c_7\sum_{\ell\geq L}e^{-c_5\ell}(1+z-K)\Big)e^{-z},
\end{eqnarray*}
where the last inequality comes from the fact that $\sum_{j\geq 0}e^{-(1-c_5)j}<\infty$. We take $L_\varepsilon=-c_8\ln \varepsilon$ so that $c_7\sum_{\ell\geq L}e^{-c_5\ell}\leq \varepsilon$ for all $L\geq L_\varepsilon$. Therefore, for any $L\geq L_\varepsilon$, $n\geq 1$ and $z\geq K\geq 0$,
\begin{equation}
\P\Big(\exists |u|=n: V(u)\leq a_n(z),\; u\not\in J_{z,K,L}^\Delta(n)\Big)\leq \Big(e^K+\varepsilon(1+z-K)\Big)e^{-z},
\end{equation}
which completes the proof. $\square$

For $b\in\z_+$, we define
\begin{equation}
\mathcal{E}_n=\mathcal{E}_n(z,b):=\{\forall k\leq n-b, \min_{u\geq w_k, |u|=n}V(u)>a_n(z)\}.
\end{equation}
We note that on the event $\mathcal{E}_n\cap\{I_n\leq a_n(z)\}$, any particle located at the leftmost position must be separated from the spine after time $n-b$.

\begin{lemma}\label{latebranching}
For any $\eta>0$ and $L>0$, there exist $K(\eta)>0$, $B(L,\eta)\geq 1$ and $N(\eta)\geq 1$ such that for any $b\geq B(L,\eta)$, $n\geq N(\eta)$ and $z\geq K\geq K(\eta)$,
\begin{equation}
\Q\Big(\mathcal{E}_n^c,\; w_n\in J^\Delta_{z,K,L}(n)\Big)\leq \eta(1+L)^2(1+z-K)n^{-3/2}.
\end{equation}
\end{lemma}

We feel free to omit the proof of Lemma \ref{latebranching} since it is just a slightly stronger version of Lemma 3.8 in \cite{elie}. It follows from the same arguments.

Let us turn to the proof of Proposition \ref{conditionedconv}. We break it up into 3 steps.

\noindent\textit{Step (I) (The conditioned convergence of $(\frac{I_n(\lfloor sn \rfloor)}{\sigma\sqrt{n}};\, 0\leq s\leq \Delta)$ for $\Delta<1$ in the non-lattice case)}

Assume that the distribution of $\mathcal{L}$ is non-lattice in this step. Recall that $a_n(z)=\frac{3}{2}\ln n-z$. The tail distribution of $I_n$ has been given in Propositions 1.3 and 4.1 of \cite{elie}, recalled as follows.
\begin{fact}[\cite{elie}]\label{tailestimation}
There exists a constant $C>0$ such that
\begin{equation}
\lim_{z\rightarrow\infty}\limsup_{n\rightarrow\infty}\Big\vert \frac{e^z}{z}\P(I_n\leq a_n(z))-C\Big\vert=0.
\end{equation}
Furthermore, for any $\varepsilon>0$, there exist $N_\varepsilon\geq 1$ and $\Lambda_\varepsilon>0$ such that for any $n\geq N_\varepsilon$ and $\Lambda_\varepsilon\leq z\leq \frac{3}{2}\ln n-\Lambda_\varepsilon$,
\begin{equation}
\Big\vert\frac{e^z}{z}\P(I_n\leq a_n(z))-C\Big\vert\leq \varepsilon.       \qquad\square
\end{equation}
\end{fact}

For any continuous functional $F: D([0,\Delta],\; \r)\rightarrow [0,1]$, it is convenient to write that
\begin{equation}
\Sigma_n(F,z):=\E\bigg[F\Big(\frac{I_n(\lfloor sn \rfloor)}{\sigma\sqrt{n}};0\leq s\leq \Delta\Big)\textbf{1}_{\{I_n\leq a_n(z)\}}\bigg].
\end{equation}
In particular, if $F\equiv 1$, $\Sigma_n(1,z)=\P(I_n\leq a_n(z))$. Thus,
\begin{equation}
\frac{\Sigma_n(F,z)}{\Sigma_n(1,z)}=\E\bigg[F\Big(\frac{I_n(\lfloor sn\rfloor)}{\sigma\sqrt{n}};0\leq s\leq \Delta\Big)\Big\vert I_n\leq a_n(z)\bigg].
\end{equation}
Let us prove the following convergence for $0<\Delta<1$,
\begin{equation}\label{latticeconv}
\lim_{z\rightarrow\infty}\limsup_{n\rightarrow\infty}\Big\vert \frac{\Sigma_n(F,z)}{\Sigma_n(1,z)}-\E[F(e_s,\;0\leq s\leq \Delta)]\Big\vert=0.
\end{equation}

\noindent\textit{Proof of (\ref{latticeconv}).} For any $n\geq 1$, $L\geq 0$ and $z\geq K\geq 0$, let
\begin{equation}
\Pi_n(F)=\Pi_n(F,z,K,L):=\E\bigg[F\Big(\frac{I_n( sn )}{\sigma\sqrt{n}};0\leq s\leq \Delta\Big)\textbf{1}_{\{ m^{(n)}\in J_{z,K,L}^\Delta(n)\}}\bigg].
\end{equation}
By Lemma \ref{typicalmin}, we obtain that for $L\geq L_\varepsilon$, $n\geq 1$ and $z\geq K\geq 0$,
\begin{equation}\label{estimationofsigma}
\Big\vert\Sigma_n(F,z)-\Pi_n(F)\Big\vert\leq \Big(e^K+\varepsilon(1+z-K)\Big)e^{-z}.
\end{equation}
Note that $m^{(n)}$ is chosen uniformly among the particles located at the leftmost position. Thus,
\begin{eqnarray}
\Pi_n(F)&=&\E\bigg[\sum_{|u|=n}\textbf{1}_{(u=m^{(n)},\;u\in J_{z,K,L}^\Delta(n))}F\Big(\frac{V(u_{\lfloor sn\rfloor})}{\sigma\sqrt{n}};0\leq s\leq \Delta\Big)\bigg]\nonumber\\
&=&\E\bigg[\frac{1}{\sum_{|u|=n}\textbf{1}_{(V(u)=I_n)}}\sum_{|u|=n}\textbf{1}_{(V(u)=I_n,\;u\in J_{z,K,L}^\Delta(n))}F\Big(\frac{V(u_{\lfloor sn\rfloor})}{\sigma\sqrt{n}};0\leq s\leq \Delta\Big)\bigg].\nonumber
\end{eqnarray}
Applying the change of measures given in (\ref{RNdensity}), it follows from Proposition \ref{changeofmeasure} that
 \begin{equation}
 \Pi_n(F) =\E_{\Q}\bigg[\frac{e^{V(w_n)}}{\sum_{|u|=n}\textbf{1}_{(V(u)=I_n)}}\textbf{1}_{(V(w_n)=I_n,\;w_n\in J_{z,K,L}^\Delta(n))}F\Big(\frac{V(w_{\lfloor sn\rfloor})}{\sigma\sqrt{n}};0\leq s\leq \Delta\Big)\bigg].
 \end{equation}
In order to estimate $\Pi_n$, we restrict ourselves to the event $\mathcal{E}_n$. Define
\begin{equation*}
\Lambda_n(F):= \E_{\Q}\bigg[\frac{e^{V(w_n)}}{\sum_{|u|=n}\textbf{1}_{(V(u)=I_n)}}\textbf{1}_{(V(w_n)=I_n,\;w_n\in J_{z,K,L}^\Delta(n))}F\Big(\frac{V(w_{\lfloor sn\rfloor})}{\sigma\sqrt{n}};0\leq s\leq \Delta\Big);\;\mathcal{E}_n\bigg].
\end{equation*}
In view of Lemma \ref{latebranching}, for any $b\geq B(L,\eta)$, $n\geq N(\eta)$ and $z\geq K\geq K(\eta)$,
\begin{eqnarray}\label{restrictiontoE}
\Big\vert\Pi_n(F)-\Lambda_n(F)\Big\vert&\leq& \E_{\Q}\Big[e^{V(w_n)};\; w_n\in J_{z,K,L}^\Delta(n),\, \mathcal{E}_n^c\Big]\\
&\leq& e^{-z}n^{-3/2}\Q\Big(\mathcal{E}_n^c,\, w_n\in J_{z,K,L}^\Delta(n)\Big)\nonumber\\
&\leq& \eta(1+L)^2(1+z-K)e^{-z}.\nonumber
\end{eqnarray}
On the event $\mathcal{E}_n\cap\{I_n\leq a_n(z)\}$, $\Lambda_n(F)$ equals
\begin{equation*}
\E_{\Q}\bigg[\frac{e^{V(w_n)}}{\sum_{u>w_{n-b},|u|=n}\textbf{1}_{(V(u)=I_n)}}\textbf{1}_{(V(w_n)=I_n,\;w_n\in J_{z,K,L}^\Delta(n))}F\Big(\frac{V(w_{\lfloor sn\rfloor})}{\sigma\sqrt{n}};0\leq s\leq \Delta\Big);\;\mathcal{E}_n\bigg].
\end{equation*}
Let, for $x\geq 0$, $L>0$, and $b\geq1$,
\begin{eqnarray}\label{latebranchingfactor}
f_{L,b}(x)&:=&\E_{\Q_x}\bigg[\frac{e^{V(w_b)-L}\textbf{1}_{\{V(w_b)=I_b\}}}{\sum_{|u|=b}\textbf{1}_{\{V(u)=I_b\}}},\min_{0\leq k\leq b}V(w_k)\geq 0, V(w_b)\leq L\bigg]\nonumber\\
&\leq& \Q_x\Big(\min_{0\leq k\leq b}V(w_k)\geq 0, V(w_b)\leq L\Big).
\end{eqnarray}
We choose $n$ large enough so that $\Delta n\leq n-b$. Thus, applying the Markov property at time $n-b$ yields that
\begin{multline}
\Lambda_n(F)=n^{3/2}e^{-z}\E_{\Q}\Big[F\Big(\frac{V(w_{\lfloor sn\rfloor})}{\sigma\sqrt{n}};0\leq s\leq \Delta\Big)f_{L,b}(V(w_{n-b})-a_n(z+L));\\
   \min_{0\leq k\leq n-b}V(w_k)\geq-z+K,\;\min_{\Delta n\leq k\leq n-b}V(w_k)\geq a_n(z+L),\;\mathcal{E}_n\Big].
\end{multline}
Let us introduce the following quantity by removing the restriction to $\mathcal{E}_n$:
\begin{multline}
\Lambda^I_n(F):= n^{3/2}e^{-z}\E_{\Q}\Big[F\Big(\frac{V(w_{\lfloor sn\rfloor})}{\sigma\sqrt{n}};0\leq s\leq \Delta\Big)f_{L,b}(V(w_{n-b})-a_n(z+L));\\
   \min_{0\leq k\leq n-b}V(w_k)\geq-z+K,\;\min_{\Delta n\leq k\leq n-b}V(w_k)\geq a_n(z+L)\Big].
\end{multline}
We immediately observe that
\begin{multline}
\Big\vert\Lambda_n(F)-\Lambda^I_n(F)\Big\vert\leq n^{3/2}e^{-z}\Q\Big(f_{L,b}(V(w_{n-b})-a_n(z+L)),\\
   \min_{0\leq k\leq n-b}V(w_k)\geq-z+K,\min_{\Delta n\leq k\leq n-b}V(w_k)\geq a_n(z+L);\; (\mathcal{E}_n)^c\Big).
\end{multline}
By (\ref{latebranchingfactor}), we check that $\Big\vert\Lambda_n(F)-\Lambda^I_n(F)\Big\vert \leq n^{3/2}e^{-z}\Q(w_n\in J_{z,K,L}^\Delta(n),\;(\mathcal{E}_n)^c)$. Applying Lemma \ref{latebranching} again implies that
\begin{equation}
\Big\vert\Lambda_n(F)-\Lambda^I_n(F)\Big\vert\leq \eta(1+L)^2(1+z-K)e^{-z}.
\end{equation}
Combining with (\ref{restrictiontoE}), we obtain that for any $b\geq B(L,\eta)$, $z\geq K\geq K(\eta)$ and $n$ large enough,
\begin{equation}\label{estimationofpi}
\Big\vert\Pi_n(F)-\Lambda^I_n(F)\Big\vert \leq 2\eta(1+L)^2(1+z-K)e^{-z}.
\end{equation}
Note that $(V(w_k);\; k\geq 1)$ is a centered random walk under $\Q$ and that it is proved in \cite{elie} that $f_{L,b}$ satisfies the conditions of Lemma \ref{keyconv}. By (I) of Lemma \ref{keyconv}, we get that
\begin{equation}
\lim_{n\rightarrow\infty}\Lambda^I_n(F)=\alpha^I_{L,b} R(z-K)e^{-z}\E[F(e_s,\; 0\leq s\leq \delta)],
\end{equation}
where $\alpha^I_{L,b}:=C_1 \int_{x\geq 0}f_{L,b}(x)R_-(x)dx\in[0,\infty)$. Thus, by (\ref{estimationofpi}), one sees that for any $b\geq B(L,\eta)$ and $z\geq K\geq K(\eta)$,
\begin{equation}
\limsup_{n\rightarrow\infty}\Big\vert\Pi_n(F)-\alpha^I_{L,b} R(z-K)e^{-z}\E[F(e_s,\; 0\leq s\leq \Delta)]\Big\vert \leq 2\eta(1+L)^2(1+z-K)e^{-z}.
\end{equation}
Going back to (\ref{estimationofsigma}), we deduce that for any $L\geq L_\varepsilon$, $b\geq B(L,\eta)$ and $z\geq K\geq K(\eta)$,
\begin{eqnarray*}
&&\limsup_{n\rightarrow\infty}\Big\vert\Sigma_n(F,z)-\alpha^I_{L,b} R(z-K)e^{-z}\E[F(e_s,\; 0\leq s\leq \Delta)]\Big\vert \\
&\leq& 2\eta(1+L)^2(1+z-K)e^{-z}+\Big(e^K+\varepsilon(1+z-K)\Big)e^{-z}.
\end{eqnarray*}
Recall that $\lim_{z\rightarrow\infty}\frac{R(z)}{z}=c_0$. We multiply each term by $\frac{e^z}{z}$, and then let $z$ go to infinity to conclude that
\begin{equation}
\limsup_{z\rightarrow\infty}\limsup_{n\rightarrow\infty}\Big\vert\frac{e^z}{z}\Sigma_n(F,z)-\alpha^I_{L,b} c_0\E[F(e_s,\; 0\leq s\leq \Delta)]\Big\vert\leq 2\eta(1+L)^2+\varepsilon.
\end{equation}
In particular, taking $F\equiv 1$ gives that
\begin{equation}
\limsup_{z\rightarrow\infty}\limsup_{n\rightarrow\infty}\Big\vert\frac{e^z}{z}\P(I_n\leq a_n(z))-\alpha^I_{L,b}c_0\Big\vert \leq 2\eta(1+L)^2+\varepsilon.
\end{equation}
It follows from Fact \ref{tailestimation} that $|C-\alpha^I_{L,b}c_0|\leq 2\eta(1+L)^2+\varepsilon$.
We thus choose $0<\varepsilon<C/10$ and $0<\eta\leq \frac{\varepsilon}{2(1+L_\varepsilon)^2}$ so that $2C>\alpha^I_{L_\varepsilon,b}c_0>C/2>0$.

Therefore, for any $\varepsilon\in(0,C/10)$, $0<\eta\leq \frac{\varepsilon}{2(1+L_\varepsilon)^2}$, $L=L_\varepsilon$ and $b\geq B(L_\varepsilon,\eta)$,
\begin{equation}
\limsup_{z\rightarrow\infty}\limsup_{n\rightarrow\infty}\bigg\vert\frac{\Sigma_n(F,z)}{\Sigma_n(1,z)}-\E[F(e_s,0\leq s\leq \Delta)]\bigg\vert\leq \frac{4\varepsilon}{C/2-2\varepsilon},
\end{equation}
which completes the proof of (\ref{latticeconv}) in the non-lattice case.

\noindent\textit{Step (II) (The conditioned convergence of $(\frac{I_n(sn)}{\sigma\sqrt{n}};\, 0\leq s\leq \Delta)$ for $\Delta<1$ in the lattice case)} Assume that the law of $\mathcal{L}$ is supported by $\alpha+\beta\z$ with span $\beta$. Recall that $a_n(0)=\alpha n+\beta\lfloor \frac{\frac{3}{2}\ln n-\alpha n}{\beta}\rfloor$ and that $a_n(z)=a_n(0)-z$. We use the same notation of Step (I). Let us prove
\begin{equation}\label{non-latticeconv}
\lim_{\beta\z\ni z\rightarrow\infty}\limsup_{n\rightarrow\infty}\Big\vert \frac{\Sigma_n(F,z)}{\Sigma_n(1,z)}-\E[F(e_s,\;0\leq s\leq \Delta)]\Big\vert=0.
\end{equation}

Suppose that $z\in\beta\z$. Whereas the arguments of Step (I), we obtain that for any $ L\geq L_\varepsilon$, $b\geq B( L,\eta)$, $z\geq K\geq K(\eta)$ and $n$ sufficiently large,
\begin{equation}\label{latticebound}
\Big\vert\Sigma_n(F,z)-\Lambda^{II}_n(F)\Big\vert \leq 2\eta(1+ L)^2(1+ z-K)e^{- z}+\Big(e^K+\varepsilon(1+z-K)\Big)e^{- z},
\end{equation}
where \begin{multline*}
\Lambda^{II}_n(F)=\Lambda^{II}(F,z,K,L,b):= e^{a_n(0)}e^{-z}\E_{\Q}\Big[F\Big(\frac{V(w_{\lfloor sn\rfloor})}{\sigma\sqrt{n}};0\leq s\leq \Delta\Big)\times\\
 f_{L,b}\big(V(w_{n-b}-a_n(z+L))\big);\; \min_{0\leq k\leq n-b}V(w_k)\geq-z+K,\;\min_{\Delta n\leq k\leq n-b}V(w_k)\geq a_n(z+L)\Big].
\end{multline*}
Under $\Q$, the distribution of $V(w_1)-V(w_0)$ is also supported by $\alpha+\beta\z$. Let $d=d(L,b):=\beta \lceil \frac{\alpha b-L}{\beta}\rceil-\alpha b+L$ and $\lambda_n:=n^{3/2}e^{-a_n(0)}$. Recall that $f_{L,b}$ is well defined in (\ref{latebranchingfactor}), it follows from (II) of Lemma \ref{keyconv} that
\begin{equation}\label{convlattice}
\lim_{n\rightarrow\infty}\lambda_n\Lambda^{II}_n(F)=\alpha_{L,b}^{II}R(z-K)e^{-z}\E[F(e_s,\; 0\leq s\leq\Delta)].
\end{equation}
where $\alpha_{L,b}^{II}:=C_1\beta\sum_{j\geq 0}f_{L,b}(\beta j+d)R_-(\beta j+d)\in[0,\infty)$. Observe that $1\leq \lambda_n\leq e^\beta$. Combining with (\ref{latticebound}), we conclude that
\begin{equation}\label{latticekey}
\limsup_{\beta\z\ni z\rightarrow\infty}\limsup_{n\rightarrow\infty}\Big\vert\frac{e^z}{z}\lambda_n\Sigma_n(F,z)-\alpha_{L,b}^{II}c_0\E[F(e_s,\;0\leq s\leq\Delta)]\Big\vert \leq e^\beta(2\eta(1+ L)^2+\varepsilon).
\end{equation}

We admit for the moment that there exist $0<c_9<c_{10}<\infty$ such that $\alpha^{II}_{L,b}\in [c_{9}, c_{10}]$ for all $L$, $b$ large enough. Then take $\varepsilon<\frac{c_{9}c_0}{4e^\beta}$, $L=L_\varepsilon$, $\eta=\frac{\varepsilon}{2(1+L_\varepsilon)^2}$ and $b\geq B(L_\varepsilon,\eta)$ so that $e^\beta(2\eta(1+ L)^2+\varepsilon)<c_{9}c_0/2\leq \alpha^{II}_{L_\varepsilon,b}c_0/2\leq 2c_{10}c_0$. Note that $\frac{\Sigma_n(F,z)}{\Sigma_n(1,z)}=\frac{\frac{e^z}{z}\lambda_n\Sigma_n(F,z)}{\frac{e^z}{z}\lambda_n\Sigma_n(1,z)}$. We thus deduce from (\ref{latticekey}) that
\begin{equation}
\limsup_{\beta\z\ni z\rightarrow\infty}\limsup_{n\rightarrow\infty} \Big\vert\frac{\Sigma_n(F,z)}{\Sigma_n(1,z)}-\E[F(e_s,\;0\leq s\leq \Delta)]\Big\vert\leq \frac{4\varepsilon}{c_{9}c_0/e^\beta-2\varepsilon},
\end{equation}
which tends to zero as $\varepsilon\downarrow 0$.

It remains to prove that $\alpha^{II}_{L,b}\in [c_{9}, c_{10}]$ for all $L$, $b$ large enough. Instead of investigating the entire system, we consider the branching random walk killed at 0. Define
\begin{equation}\label{minkilling}
I^{kill}_n:=\inf\{V(u): |u|=n, V(u_k)\geq 0,\; \forall 0\leq k\leq n\},
\end{equation}
and we get the following fact from Corollary 3.4 and Lemma 3.6 of \cite{elie}.
\begin{fact}[\cite{elie}]\label{upper}
There exists a constant $c_{11}>0$ such that for any $n\geq 1$ and $x,\; z\geq 0$,
\begin{equation}\label{upperboundkilling}
\P_x(I_n^{kill}\leq a_n(z))\leq c_{11}(1+x)e^{-x-z}.
\end{equation}
Moreover, there exists $c_{12}>0$ such that for any $n\geq 1$ and $z\in[0, a_n(1)]$,
\begin{equation}\label{lowerboundkilling}
\P(I^{kill}_n\leq a_n(z))\geq c_{12} e^{-z}.
\end{equation}
\end{fact}
Even though Fact \ref{upper} is proved in \cite{elie} under the assumption that the distribution of $\mathcal{L}$ is non-lattice, the lattice case is actually recovered from that proof.

{\color{black}Analogically, let $m^{kill,(n)}$ be the particle chosen uniformly in the set $\{u: |u|=n,\;V(u)=I_n^{kill},\; \min_{0\leq k\leq n}V(u_k)\geq 0\}$. Moreover, let $\Sigma^{kill}_n(1,z):=\P\Big[I^{kill}_n\leq a_n(z)\Big]$ and $\Pi_n^{kill}(1,z,z,L):=\P\Big[I_n^{kill}\leq a_n(z),\; m^{kill,(n)}\in J^\Delta_{z,z,L}(n)\Big]$. By (\ref{killingupp}) again, we check that for all $L\geq L_\varepsilon$,
\begin{eqnarray}\label{eq:killing1}
&&\Big\vert\Sigma_n^{kill}(1,z)-\Pi_n^{kill}(1,z,z,L)\Big\vert \\
& \leq &\P\Big[\exists |u|=n: V(u)\leq a_n(z);\; \min_{0\leq k\leq n}V(u_k)\geq 0;\; \min_{\Delta n\leq k\leq n}V(u_k)\leq a_n(z+L)\Big]\nonumber\\
&\leq &\varepsilon e^{-z}.\nonumber
\end{eqnarray}
Recounting the arguments of Step (1), one sees that for any $L\geq L_\varepsilon$, $b\geq B(L,\eta)$, $z\geq K(\eta)$ and $n$ sufficiently large,
\begin{equation}\label{latticeboundkilling}
\Big\vert \Pi^{kill}_n(1,z,z,L)-\Lambda^{kill}_n\Big\vert \leq 2\eta(1+ L)^2e^{-z} ,
\end{equation}
where
\begin{equation}
\Lambda^{kill}_n:= \E_{\Q}\Big[f^{kill}(V(w_{n-b})); \min_{0\leq k\leq n-b}V(w_k)\geq 0,\min_{\Delta n\leq k\leq n-b}V(w_k)\geq a_n(z+L)\Big],
\end{equation}
with $f^{kill}(x):=\E_{\Q_x}\Big[\frac{e^{V(w_b)}\textbf{1}_{\{V(w_b)=I_b^{kill}\}}}{\sum_{|u|=b}\textbf{1}_{\{V(u)=I_b^{kill},\;\min_{0\leq j\leq b}V(u_j)\geq 0\}}};\;\min_{0\leq k\leq b}V(w_k)\geq a_n(z+L), V(w_b)\leq a_n(z)\Big]$. For $\varepsilon>0$ and $n$ sufficiently large, it has been proved in \cite{elie} that
\begin{equation}\label{killestim}
\Big\vert e^z \Lambda^{II}_n(1,z,z,L,b)-\Lambda^{kill}_n\Big\vert\leq \varepsilon.
\end{equation}
Recalling the convergence (\ref{convlattice}) with $K=z$ and $F\equiv1$, we deduce from (\ref{eq:killing1}), (\ref{latticeboundkilling}) and (\ref{killestim}) that for any $L\geq L_\varepsilon$, $b\geq B(L,\eta)$ and $z\geq K(\eta)$,
\begin{equation}
\limsup_{n\rightarrow\infty}\Big\vert\lambda_n \Sigma_n^{kill}(1,z)-\alpha_{L,b}^{II}e^{-z}\Big\vert\leq e^{\beta}\Big(2\eta(1+ L)^2+2\varepsilon\Big)e^{-z},
\end{equation}
since $R(0)=1$ and $1\leq \lambda_n\leq e^\beta$. Fact \ref{upper} implies that $c_{12}\leq e^z\lambda_n\P(I_n^{kill}\leq a_n(z))\leq c_{11} e^\beta$. Hence, we obtain that
\begin{equation}
c_{12}-e^\beta\Big(2\eta(1+ L)^2 +2 \varepsilon\Big)\leq \alpha^{II}_{L,b} \leq e^{\beta}c_{11}+e^\beta\Big(2\eta(1+ L)^2 +2 \varepsilon\Big).
\end{equation}
Let $c_{10}:=c_{11}e^{\beta}+c_{12}$ and $c_{9}:=3c_{12}/4>0$. For any $\varepsilon<e^{-\beta}c_{12}/12$, we take $L=L_\varepsilon$ and $\eta\leq \varepsilon/2(1+L_\varepsilon)^2$. Then $c_{10}>\alpha_{L,b}^{II}\geq c_{9}>0$ for $b\geq B(L_\varepsilon,\eta)$. This completes the second step.
}

\noindent\textit{Step (III)(The tightness)} Actually, it suffices to prove the following proposition.
\begin{proposition}\label{tightness at end}
For any $\eta>0$,
\begin{equation}
\lim_{\delta\rightarrow0}\limsup_{z\rightarrow\infty}\limsup_{n\rightarrow\infty}\P\Big(\sup_{0\leq k\leq \delta n}|I_n(n-k)-I_n|\geq \eta \sigma\sqrt{n}\Big\vert I_n\leq a_n(z) \Big)=0.
\end{equation}
\end{proposition}

The first two steps allow us to obtain the following fact whether the distribution is lattice or non-lattice.

\begin{fact}\label{taillaw}
There exist constants $c_{13},c_{14}\in(0,\infty)$ such that
\begin{equation}
c_{13}\leq\liminf_{z\rightarrow\infty}\liminf_{n\rightarrow\infty}\frac{e^{z}}{z}\P(I_n\leq a_n(z))\leq \limsup_{z\rightarrow\infty}\limsup_{n\rightarrow\infty}\frac{e^{z}}{z}\P(I_n\leq a_n(z))\leq c_{14}.
\end{equation}
\end{fact}

\noindent\textit{Proof of Proposition \ref{tightness at end}.}
First, we observe that for any $M\geq 1$ and $\delta\in(0,1/2)$,
\begin{eqnarray}\label{simpleob}
&&\P\Big(\sup_{0\leq k\leq \delta n}|I_n(n-k)-I_n|\geq \delta\sigma \sqrt{n},\; I_n\leq a_n(z)\Big)\nonumber\\
&\leq& \P\Big(m^{(n)}_n\not\in J_{z,0,L}^{1/2}(n),\, I_n\leq a_n(z)\Big)+\P\Big(I_n(n-\lfloor\delta n\rfloor)\geq M\sigma\sqrt{\delta n},\, I_n\leq a_n(z)\Big)+\chi(\delta,z,n)\nonumber.
\end{eqnarray}
where $\chi(\delta,z,n):= \P\Big(\, m_n^{(n)}\in J_{z,0,L}^{1/2}(n),\, I_n(n-\lfloor\delta n\rfloor)\leq M\sigma\sqrt{\delta n},\,\sup_{0\leq k\leq \delta n}|I_n(n-k)-I_n|\geq \eta\sigma \sqrt{n}\Big)$.

It follows from Lemma \ref{typicalmin} that for any $\varepsilon>0$, if $L\geq L_\varepsilon$, $n\geq 1$ and $z\geq 0$,
\begin{equation}
\P\Big(m_n^{(n)}\not\in J_{z,0,L}^{1/2}(n),\, I_n\leq a_n(z)\Big)\leq (1+\varepsilon (1+z))e^{-z}.
\end{equation}
Then dividing each term of (\ref{simpleob}) by $\P(I_n\leq a_n(z))$ yields that
\begin{eqnarray}
&&\P\Big(\sup_{0\leq k\leq \delta n}|I_n(n-k)-I_n|\geq \eta\sigma \sqrt{n}\Big\vert I_n\leq a_n(z)\Big)\\
&\leq& \frac{(1+\varepsilon (1+z))e^{-z}}{\P(I_n\leq a_n(z))}+\P\Big(I_n(n-\lfloor\delta n\rfloor)\geq M\sigma\sqrt{\delta n}\Big\vert I_n\leq a_n(z)\Big)+\frac{\chi(\delta,z,n)}{\P(I_n\leq a_n(z))}.\nonumber
\end{eqnarray}
On the one hand, by Fact \ref{taillaw},
\begin{equation}\label{term1}
\limsup_{z\rightarrow\infty}\limsup_{n\rightarrow\infty} \frac{(1+\varepsilon (1+z))e^{-z}}{\P(I_n\leq a_n(z))}\leq \frac{\varepsilon}{c_{13}}.
\end{equation}
On the other hand, Steps (I) and (II) tell us that for any $1>\delta>0$ and $M\geq 1$,
\begin{equation}\label{upperbound}
\limsup_{z\rightarrow\infty}\limsup_{n\rightarrow\infty}\P\Big[I_n(n-\lfloor\delta n\rfloor)\geq M\sigma\sqrt{\delta n}\Big\vert I_n\leq a_n(z) \Big]=\P[e_{1-\delta}\geq M\sqrt{\delta}],
\end{equation}
which, by Chebyshev's inequality, is bounded by $\frac{\E[e_{1-\delta}]}{ M\sqrt{\delta}}= \frac{4\sqrt{1-\delta}}{M\sqrt{2\pi}}$. Consequently,
\begin{eqnarray}\label{detailestimate}
&&\limsup_{z\rightarrow\infty}\limsup_{n\rightarrow\infty}\P\Big(\sup_{0\leq k\leq \delta n}|I_n(n-k)-I_n|\geq \eta\sigma \sqrt{n}\Big\vert I_n\leq a_n(z)\Big)\\
&\leq& \frac{\varepsilon}{c_{13}}+\frac{2}{M}+\limsup_{z\rightarrow\infty}\limsup_{n\rightarrow\infty}\frac{\chi(\delta,z,n)}{\P(I_n\leq a_n(z))}.\nonumber
\end{eqnarray}

Let us estimate $\chi(\delta,z,n)$. One sees that
\begin{eqnarray*}
\chi(\delta,z,n)&\leq &\E\Big[\sum_{|u|=n}\textbf{1}_{\{u\in J_{z,L}^{1/2}(n);\;\sup_{0\leq k\leq\delta n}|V(u_{n-k})-V(u)|\geq \eta\sigma\sqrt{n};\; V(u_{n-\lfloor \delta n\rfloor})\leq M\sigma\sqrt{\delta n}\}}\Big].
\end{eqnarray*}
By Lemma \ref{many-to-one}, it becomes that
\begin{eqnarray*}
\chi(\delta,z,n)&\leq&\E\Big[e^{S_n}; S_n\leq a_n(z), \underline{S}_n\geq-z, \underline{S}_{[n/2, n]}\geq a_n(z+L),\\
 &&\qquad S_{n-\lfloor \delta n\rfloor}\leq M\sigma\sqrt{\delta n},\sup_{0\leq k\leq \delta n}|S_{n-k}-S_n|\geq \eta\sigma\sqrt{n}\Big]\\
&\leq & n^{3/2}e^{-z}\Upsilon(\delta,z,n),
\end{eqnarray*}
where $\Upsilon(\delta,z,n):=\P\Big(S_n\leq a_n(z),\;\underline{S}_n\geq-z,\;  \underline{S}_{[n/2, n]}\geq a_n(z+L),\;  S_{n-\lfloor \delta n\rfloor}\leq M\sigma\sqrt{\delta n},\\ \sup_{0\leq k\leq \delta n}|S_{n-k}-S_n|\geq \eta\sigma\sqrt{n},\,  S_{n-\lfloor \delta n\rfloor}\leq M\sigma\sqrt{\delta n}\Big)$.

Reversing time yields that
\begin{multline}
\Upsilon(\delta,z,n)\leq\P\Big(\underline{S}^-_n\geq-a_n(0),\, \underline{S}^-_{ n/2}\geq -L,\, -S_n\in [-a_n(z),-a_n(z+L)],\\
 \sup_{0\leq k\leq \delta n}|-S_k|\geq \eta\sigma\sqrt{n}, -S_{\lfloor \delta n\rfloor}\leq M\sigma\sqrt{\delta n}-a_n(z+L)\Big).
\end{multline}
Applying the Markov property at time $\lfloor \delta n\rfloor$, we obtain that
\begin{equation}\label{reversingtime}
\Upsilon(\delta,z,n)=\E\Big[\Theta(-S_{\lfloor\delta n\rfloor});\;\underline{S}^-_{ \delta n}\geq -L, \sup_{0\leq k\leq \delta n}|-S_k|\geq \eta\sigma\sqrt{n}\Big],
\end{equation}
where $\Theta(x):=\textbf{1}_{\{x\leq M\sigma\sqrt{\delta n}-a_n(z+L)\}}\P_x\Big(\underline{S}^-_{ (1/2-\delta)n}\geq -L, \underline{S}^-_{ (1-\delta )n}\geq -a_n(0), -S_{n-\lfloor\delta n\rfloor}\in[-a_n(z),\\
-a_n(z+L)]\Big)$.
Reversing time again implies that
\begin{multline*}
\Theta(x)\leq\textbf{1}_{\{x\leq M\sigma\sqrt{\delta n}\}}\P\Big(\underline{S}_{(1-\delta)n}\geq -z-L,\\
\underline{S}_{[n/2, (1-\delta)n]}\geq a_n(z+2L) ,S_{n-\lfloor\delta n\rfloor}\in [x+a_n(z+L), x+a_n(z)]\Big).
\end{multline*}
By (\ref{elementary}), $\Theta(x)\leq c_{15}(1+z+L)(1+L)(1+M\sigma\sqrt{\delta n}+2L)n^{-3/2}$. Plugging it into (\ref{reversingtime}) and taking $n$ large enough so that $1+2L<\eta\sigma\sqrt{\delta n}$, we get that
\begin{equation*}
\Upsilon(\delta,z,n)\leq c_{15}(1+z)(1+L)^2n^{-3/2}(M+\eta)\sigma\sqrt{\delta n}\E\Big[\underline{S}^-_{\delta n}\geq -L, \sup_{0\leq k\leq \delta n}|-S_k|\geq \eta\sigma\sqrt{n}\Big].
\end{equation*}
Recall that $\chi(\delta,z,n)\leq e^{-z}n^{3/2}\Upsilon(\delta,z,n)$. We check that
\begin{multline}
\chi(\delta,z,n)\leq c_{15}e^{-z}(1+z)(1+L)^2(M+\eta)\sigma\\
\times\E_L\Big[ \sup_{0\leq k\leq \delta n}(-S_k)\geq \eta\sigma\sqrt{n}\Big\vert \underline{S}^-_{\delta n}\geq 0\Big]\Big(\sqrt{\delta n}\P_L\Big[\underline{S}^-_{\delta n}\geq 0\Big]\Big).
\end{multline}
On the one hand, by Theorem 1.1 of \cite{caravenna-loic}, $\E_L\Big[ \sup_{0\leq k\leq \delta n}(-S_k)\geq \eta\sigma\sqrt{n}\Big\vert \underline{S}^-_{\delta n}\geq 0\Big]$ converges to $\P(\sup_{0\leq s\leq 1}\mathcal{M}_s\geq \eta/\sqrt{\delta})$ as $n\rightarrow\infty$. On the other hand, (\ref{probstaypositive}) shows that $\sqrt{\delta n}\P_L\Big[\underline{S}^-_{\delta n}\geq 0\Big]$ converges to $C_-R_-(L)$ as $n\rightarrow\infty$.  Therefore,
\begin{multline*}\label{mainfactor}
\limsup_{n\rightarrow\infty}\chi(\delta,z,n)\leq c_{15}e^{-z}(1+z)(1+L)^2(M+\eta)\sigma C_-R_-(L)\times\P(\sup_{0\leq s\leq 1}\mathcal{M}_s\geq \eta/\sqrt{\delta}).
\end{multline*}
Going back to (\ref{detailestimate}) and letting $z\rightarrow\infty$, we deduce from Fact \ref{taillaw} that
\begin{multline}
\limsup_{z\rightarrow\infty}\limsup_{n\rightarrow\infty}\P\Big(\sup_{0\leq k\leq \delta n}|I_n(n-k)-I_n|\geq \eta\sigma \sqrt{n}\Big\vert I_n\leq a_n(z)\Big)\\
\leq \frac{\varepsilon}{c_{13}}+\frac{2}{M}+\frac{c_{15}(1+L)^2(M+\eta)\sigma C_-R_-(L)\times\P(\sup_{0\leq s\leq 1}\mathcal{M}_s\geq \eta/\sqrt{\delta})}{c_{13}}.
\end{multline}
Notice that $\P(\sup_{0\leq s\leq 1}\mathcal{M}_s\geq \eta/\sqrt{\delta})$ decreases to 0 as $\delta\downarrow0$. Take $M\geq 2/\varepsilon$. We conclude that for any $0<\varepsilon<c_{13}$,
\begin{equation}
\limsup_{\delta\rightarrow0}\limsup_{z\rightarrow\infty}\limsup_{n\rightarrow\infty}\P\Big(\sup_{0\leq k\leq \delta n}|I_n(n-k)-I_n|\geq \eta\sigma \sqrt{n}\Big\vert I_n\leq a_n(z)\Big)\leq \frac{\varepsilon}{c_{13}}+\varepsilon,
\end{equation}
which completes the proof of Proposition \ref{tightness at end}. And Proposition \ref{conditionedconv} is thus proved.  $\square$

\section{Proof of Theorem \ref{mainthm}}
Let us prove the main theorem now. It suffices to prove that for any continuous functional $F: D([0,1],\r)\rightarrow[0,1]$, we have
\begin{equation}\label{trueconv}
\lim_{n\rightarrow\infty}\bigg\vert\E\Big[F\Big(\frac{I_n(\lfloor sn \rfloor)}{\sigma\sqrt{n}};0\leq s\leq 1\Big) \Big]-\E\Big[F(e_s,\;0\leq s\leq1)\Big]\bigg\vert=0.
\end{equation}

\textit{Proof of (\ref{trueconv}).} Define for $A\geq 0$,
\begin{equation}
\mathcal{Z}[A]:=\{u\in\mathbb{T}: V(u)\geq A>\max_{k<|u|}V(u_k)\}.
\end{equation}
For any particle $u\in \mathcal{Z}[A]$, there is a subtree rooted at $u$. If $|u|\leq n$, let
\begin{equation*}
I_{n}(u):=\min_{v\geq u, |v|=n}V(v).
\end{equation*}
Moreover, assume $m_n^u$ is the particle uniformly chosen in the set $\{|v|=n: v\geq u, V(v)=I_n(u)\}$. Similarly, we write $[\![\varnothing,m_n^u]\!]:=\{\varnothing=:m_0^u, m_1^u,\cdots, m_n^u\}$. The trajectory leading to $m_n^u$ is denoted by $\{V(m_k^u);0\leq k\leq n\}$. Let $\omega_A$ be the particle uniformly chosen in $\{u\in\mathcal{Z}[A]: |u|\leq n,\, I_{n}(u)=I_n\}$.

Let $\mathcal{Y}_A:=\{\max_{u\in\mathcal{Z}[A]}|u|\leq M,\, \max_{u\in\mathcal{Z}[A]}V(u)\leq M\}$. Then for any $\varepsilon>0$, there exist $M:=M(A,\varepsilon)$ large enough such that $\P(\mathcal{Y}_A^c)\leq \varepsilon$. It follows that
\begin{eqnarray}\label{rearrange}
&&\Big\vert\E\Big[F\Big(\frac{I_n(\lfloor sn\rfloor)}{\sigma\sqrt{n}};0\leq s\leq 1\Big)\Big]-\E\Big[F\Big(\frac{I_n(\lfloor sn\rfloor)}{\sigma\sqrt{n}};0\leq s\leq 1\Big); \mathcal{Y}_A, |I_n- a_n(0)|\leq A/2\Big]\Big\vert\\
&\leq&\varepsilon+\P[|I_n- a_n(0)|\geq A/2].\nonumber
\end{eqnarray}
We then check that for $n\geq M$,
\begin{eqnarray}\label{rearrangement}
&&\E\Big[F\Big(\frac{I_n(\lfloor sn\rfloor)}{\sigma\sqrt{n}};0\leq s\leq 1\Big); \mathcal{Y}_A, |I_n-a_n(0)|\leq A/2\Big]\\
&=&\E\Big[\sum_{u\in\mathcal{Z}[A]}\textbf{1}_{(u=\omega_A)}F\Big(\frac{V(m^u_{\lfloor sn\rfloor})}{\sigma\sqrt{n}};0\leq s\leq 1\Big); \mathcal{Y}_A, |I_n- a_n(0)|\leq A/2\Big].\nonumber
\end{eqnarray}
Define another trajectory $\{\tilde{V}(m_k^u);0\leq k\leq n\}$ as follows.
\begin{equation}
\tilde{V}(m_k^u):=\begin{cases}
 V(u) & \text{if } k<|u|;\\
 V(m_k^u) & \text{if }|u|\leq k\leq n.
\end{cases}
\end{equation}
It follows that
\begin{eqnarray}\label{rearrangement1}
&&\E\Big[F\Big(\frac{I_n(\lfloor sn\rfloor)}{\sigma\sqrt{n}};0\leq s\leq 1\Big); \mathcal{Y}_A, |I_n- a_n(0)|\leq A/2\Big]\\
&=&\E\Big[\sum_{u\in\mathcal{Z}[A]}\textbf{1}_{(u=\omega_A)}F\Big(\frac{\tilde{V}(m^u_{\lfloor sn\rfloor})}{\sigma\sqrt{n}};0\leq s\leq 1\Big); \mathcal{Y}_A, |I_n- a_n(0)|\leq A/2\Big]+o_n(1),\nonumber
\end{eqnarray}
where $o_n(1)\rightarrow 0$ as $n$ goes to infinity.

Define the sigma-field $\mathcal{G}_A:=\sigma\{(u,V(u),I_n(u)); u\in\mathcal{Z}[A]\}$.
Note that on $\mathcal{Y}_A$, $I_n=\min_{u\in\mathcal{Z}[A]}I_n(u)$ as long as $n\geq M$.
One sees that $\mathcal{Y}_A\cap\{|I_n-a_n(0)|\leq A/2\}$ is $\mathcal{G}_A$-measurable for all $n$ large enough.
Thus,
\begin{eqnarray}
&&\E\Big[\sum_{u\in\mathcal{Z}[A]}\textbf{1}_{(u=\omega_A)}F\Big(\frac{\tilde{V}(m^u_{\lfloor sn\rfloor})}{\sigma\sqrt{n}};0\leq s\leq 1\Big); \mathcal{Y}_A, |I_n- a_n(0)|\leq A/2\Big]\\
&=&\E\Big[\sum_{u\in\mathcal{Z}[A]}\textbf{1}_{(u=\omega_A)}\E\Big[F\Big(\frac{\tilde{V}(m^u_{\lfloor sn\rfloor})}{\sigma\sqrt{n}};0\leq s\leq 1\Big)\Big\vert\mathcal{G}_A,\, u=\omega_A\Big]; \mathcal{Y}_A, |I_n- a_n(0)|\leq A/2\Big].\nonumber
\end{eqnarray}
Further, we notice by the branching property that conditioned on $\{(u,V(u));u\in \mathcal{Z}[A]\}$, the subtrees generated by $u\in\mathcal{Z}[A]$ are independent copies of the original one, started from $V(u)$, respectively. Therefore, given $\mathcal{Y}_A\cap\{ |I_n- a_n(0)|\leq A/2\}$,
\begin{eqnarray*}
&&\textbf{1}_{(u=\omega_A)}\E\Big[F\Big(\frac{\tilde{V}(m^u_{\lfloor sn\rfloor})}{\sigma\sqrt{n}};0\leq s\leq 1\Big)\Big\vert\mathcal{G}_A,\, u=\omega_A\Big]\\
&=&\textbf{1}_{(u=\omega_A)}\E\Big[F\Big(\frac{I({\lfloor s(n-|u|)\rfloor})}{\sigma\sqrt{n-|u|}};0\leq s\leq 1\Big)\Big\vert I_{n-|u|}\leq a_n(-r_u)\Big]+o_n(1),
\end{eqnarray*}
where $r_u:=\min\{\min_{v\in\mathcal{Z}[A]\setminus\{u\}}I_n(v)-a_n(0), \, A/2\}-V(u)$ is independent of $I_{n-|u|}$. Thus, (\ref{rearrangement1}) becomes that
\begin{eqnarray}\label{rearrangement2}
&&\E\Big[F\Big(\frac{I_n(\lfloor sn\rfloor)}{\sigma\sqrt{n}};0\leq s\leq 1\Big); \mathcal{Y}_A, |I_n- a_n(0)|\leq A/2\Big]\\
&=&\E\Big[\sum_{u\in\mathcal{Z}[A]}\textbf{1}_{(u=\omega_A)}\E\Big[F\Big(\frac{I({\lfloor s(n-|u|)\rfloor})}{\sigma\sqrt{n-|u|}};0\leq s\leq 1\Big)\Big\vert I_{n-|u|}\leq a_n(-r_u)\Big];\nonumber \\
&&\qquad \mathcal{Y}_A, |I_n- a_n(0)|\leq A/2\Big]+o_n(1).\nonumber
\end{eqnarray}
The event  $\mathcal{Y}_A\cap\{|I_n-a_n(0)|\leq A/2\}$ ensures that $A/2+M\geq -r_u\geq A/2$. The conditioned convergence has been given in Proposition \ref{conditionedconv}. We need a slightly stronger version here.

According to Proposition \ref{conditionedconv}, for any $\varepsilon>0$, there exists $z_\varepsilon>0$ such that for all $z\geq z_\varepsilon$,
\begin{multline}
\limsup_{n\rightarrow\infty}\Big\vert\E\Big[F\Big(\frac{I_n(\lfloor sn\rfloor)}{\sigma\sqrt{n}};0\leq s\leq 1\Big)\Big\vert I_n\leq a_n(z) \Big]-\E[F(e_s,0\leq s\leq 1)] \Big\vert<\varepsilon.
\end{multline}
Thus, for any $z\geq z_\varepsilon$, there exists $N_z\geq 1$ such that for any $n\geq N_z$,
\begin{equation}\label{conditionedestimate}
\Big\vert\E\Big[F\Big(\frac{I_n(\lfloor sn\rfloor)}{\sigma\sqrt{n}};0\leq s\leq 1\Big)\Big\vert I_n\leq a_n(z) \Big]-\E[F(e_s,0\leq s\leq 1)] \Big\vert<2\varepsilon.
\end{equation}
Take $A=2z_\varepsilon$ and $K=M$. We say that for $n$ sufficiently large,
\begin{equation}\label{uniformconv}
\sup_{z\in[z_\varepsilon, z_\varepsilon+K]}\Big\vert\E\Big[F\Big(\frac{I({\lfloor s(n)\rfloor})}{\sigma\sqrt{n}};0\leq s\leq 1\Big)\Big\vert I_{n}\leq a_n(z)\Big]-\E[F(e_s, 0\leq s\leq 1)]\Big\vert\leq 3\varepsilon.
\end{equation}
In the lattice case, (\ref{uniformconv}) follows immediately. We only need to prove it in the non-lattice case.

Recall that $\Sigma_n(F,z)=\E\Big[F\Big(\frac{I_n(\lfloor sn\rfloor)}{\sigma\sqrt{n}};0\leq s\leq 1\Big);\; I_n\leq a_n(z) \Big]$ with $0\leq F\leq 1$.
Then, for any $\ell>0$ and $z\geq0$,
\begin{eqnarray}\label{smalldiff}
&&\Big\vert\frac{\Sigma_n(F,z)}{\Sigma_n(1,z)}-\frac{\Sigma_n(F,z+\ell)}{\Sigma_n(1,z+\ell)}\Big\vert\\
&\leq& \Big\vert\frac{\Sigma_n(F,z)-\Sigma_n(F,z+\ell)}{\Sigma_n(1,z)}\Big\vert+\Big\vert
\frac{\Sigma_n(F,z+\ell)}{\Sigma_n(1,z)}-\frac{\Sigma_n(F,z+\ell)}{\Sigma_n(1,z+\ell)}\Big\vert\nonumber\\
&=&\frac{1}{\Sigma_n(1,z)}\Big(\Big\vert\Sigma_n(F,z)-\Sigma_n(F,z+\ell)\Big\vert+
\frac{\Sigma_n(F,z+\ell)}{\Sigma_n(1,z+\ell)}\Big\vert\Sigma_n(1,z+\ell)-\Sigma_n(1,z)\Big\vert\Big).\nonumber
\end{eqnarray}
Since $0\leq F\leq 1$, the two following inequalities
\begin{eqnarray*}
\Big\vert\Sigma_n(F,z)-\Sigma_n(F,z+\ell)\Big\vert&=&\E\Big[F\Big(\frac{I_n(\lfloor sn\rfloor)}{\sigma\sqrt{n}};0\leq s\leq 1\Big);\;a_n(z+\ell)<I_n\leq a_n(z)\Big]\\
&\leq & \P(a_n(z+\ell)<I_n\leq a_n(z)),
\end{eqnarray*}
and $\frac{\Sigma_n(F,z+\ell)}{\Sigma_n(1,z+\ell)}\leq 1$ hold. Note also that $\vert\Sigma_n(1,z+\ell)-\Sigma_n(1,z)\vert=\P(a_n(z+\ell)<I_n\leq a_n(z))$. It follows that
\begin{eqnarray}\label{distance}
\Big\vert\frac{\Sigma_n(F,z)}{\Sigma_n(1,z)}-\frac{\Sigma_n(F,z+\ell)}{\Sigma_n(1,z+\ell)}\Big\vert&\leq&2\frac{\P(a_n(z+\ell)<I_n\leq a_n(z))}{\P(I_n\leq a_n(z))}\\
&=&2-2\frac{\P(I_n\leq a_n(z+\ell))}{\P(I_n\leq a_n(z))}.\nonumber
\end{eqnarray}
In view of Fact \ref{tailestimation}, we take $\frac{3}{2}\ln n-\Lambda_{\varepsilon^\prime}\geq \ell+z>z\geq \Lambda_{\varepsilon^\prime}$ so that for any $n\geq N_{\varepsilon^\prime}$,
\begin{equation}
\frac{\P(I_n\leq a_n(z+\ell))}{\P(I_n\leq a_n(z))}\geq \frac{(C-\varepsilon^\prime)(z+\ell)e^{-z-\ell}}{(C+\varepsilon^\prime)ze^{-z}}\geq\frac{C-\varepsilon^\prime}{C+\varepsilon^\prime}e^{-\ell}.
\end{equation}

For $\varepsilon^\prime=C\varepsilon/8>0$, we choose $\zeta=\frac{\varepsilon}{4}$ so that $\frac{C-\varepsilon^\prime}{C+\varepsilon^\prime}e^{-\zeta}\geq 1-\frac{\varepsilon}{2}$. As a consequence, for any $ \Lambda_{\varepsilon^\prime}\leq z\leq \frac{3}{2}\ln n-\Lambda_{\varepsilon^\prime}-\zeta$, $0\leq \ell\leq \zeta$ and $n\geq N_{\varepsilon^\prime}$,
\begin{equation}\label{distance}
\Big\vert\frac{\Sigma_n(F,z)}{\Sigma_n(1,z)}-\frac{\Sigma_n(F,z+\ell)}{\Sigma_n(1,z+\ell)}\Big\vert\leq 2\Big(1-\frac{C-\varepsilon^\prime}{C+\varepsilon^\prime}e^{-\ell}\Big)\leq \varepsilon.
\end{equation}

For $\varepsilon>0$, $z_\varepsilon$ can be chosen so that $[z_\varepsilon, z_\varepsilon+K]\subset [\Lambda_{\varepsilon^\prime}, \frac{3}{2}\ln n-\Lambda_{\varepsilon^\prime}]$ for $n\geq e^KN_{\varepsilon^\prime}$. For any integer $0\leq j\leq \lceil K/\zeta\rceil$, let $z_j:=z_\varepsilon+j\zeta$. Then $[z_\varepsilon,z_\varepsilon+K]\subset\cup_{0\leq j\leq \lceil K/\zeta\rceil}[z_j,z_{j+1}]$. Take $N^\prime_\varepsilon=\max_{0\leq j\leq \lceil K/\zeta\rceil}\{N_{z_j}, e^KN_{\varepsilon^\prime}\}$. By (\ref{conditionedestimate}) and (\ref{distance}), we conclude that for any $n\geq N_\varepsilon^\prime$,
\begin{eqnarray*}
&&\sup_{z\in[z_\varepsilon,z_\varepsilon+K]}\Big\vert\E\Big[F\Big(\frac{I_n(\lfloor sn\rfloor)}{\sigma\sqrt{n}};0\leq s\leq 1\Big)\Big\vert I_n\leq a_n(z) \Big]-\E[F(e_s,0\leq s\leq 1)] \Big\vert\\
&\leq & \sup_{0\leq j\leq \lceil K/\zeta\rceil}\Big\vert \frac{\Sigma_n(F,z_j)}{\Sigma_n(1, z_j)}-\E[F(e_s,\;0\leq s\leq 1)]\Big\vert+ \sup_{0\leq j< \lceil K/\zeta\rceil}\sup_{z_j\leq z\leq z_{j+1}}\Big\vert \frac{\Sigma_n(F,z)}{\Sigma_n(1,z)}-\frac{\Sigma_n(F,z_j)}{\Sigma_n(1,z_j)}\Big\vert\\
&\leq&3\varepsilon.
\end{eqnarray*}
We continue to prove the main theorem. Since $\sum_{u\in\mathcal{Z}[A]}\textbf{1}_{(u=\omega_A)}=1$, we deduce from (\ref{rearrangement2}) and (\ref{uniformconv}) that for $n$ sufficiently large,
\begin{eqnarray*}
&&\Big\vert\E\Big[F\Big(\frac{I_n(\lfloor sn\rfloor)}{\sigma\sqrt{n}};0\leq s\leq 1\Big); \mathcal{Y}_A, |I_n- a_n(0)|\leq A/2\Big]-\E[F(e_s,\;0\leq s\leq 1)]\Big\vert\\
&\leq&3\varepsilon \P(\mathcal{Y}_A; |I_n-a_n(0)|\leq A/2)+ o_n(1)+\P(\mathcal{Y}_A^c)+\P(|I_n-a_n(0)|\geq A/2)\\
&\leq &4\varepsilon +o_n(1)+\P(|I_n-a_n(0)|\geq A/2).
\end{eqnarray*}
Going back to (\ref{rearrange}), we conclude that for $n$ large enough,
\begin{equation}
\Big\vert\E\Big[F\Big(\frac{I_n(\lfloor sn\rfloor)}{\sigma\sqrt{n}};0\leq s\leq 1\Big)\Big]-\E[F(e_s,\;0\leq s\leq 1)]\Big\vert\leq5\varepsilon+2\P(|I_n- a_n(0)|\geq A/2)+o_n(1).\nonumber
\end{equation}
Let $n$ go to infinity and then make $\varepsilon\downarrow0$. Therefore,
\begin{eqnarray}
&&\limsup_{n\rightarrow\infty}\Big\vert\E\Big[F\Big(\frac{I_n(\lfloor sn\rfloor)}{\sigma\sqrt{n}};0\leq s\leq 1\Big)\Big]-\E[F(e_s,\;0\leq s\leq 1)]\Big\vert\\
&\leq &\limsup_{z\rightarrow\infty}\limsup_{n\rightarrow\infty}2\P(|I_n- a_n(0)|\geq z).\nonumber
\end{eqnarray}
It remains to show that $\limsup_{z\rightarrow\infty}\limsup_{n\rightarrow\infty}\P(|I_n- a_n(0)|\geq z)=0$. Because of Fact (\ref{taillaw}), it suffices to prove that
\begin{equation}
\limsup_{z\rightarrow\infty}\limsup_{n\rightarrow\infty}\P(I_n\geq a_n(0)+ z)=0.
\end{equation}
In the non-lattice case, Theorem 1.1 of \cite{elie} implies it directly. In the lattice case, we see that for $n$ large enough,
\begin{equation}\label{latticetight}
\P(I_n\geq a_n(0)+ z)\leq \E\Big[\prod_{u\in\mathcal{Z}[A]}(1-\Phi_u(z,n)); \mathcal{Y}_A\Big]+\varepsilon,
\end{equation}
with $\Phi_u(z, n):=\P(I_{n-|u|}\leq a_n(V(u)-z))$. Take $A=2z$ here. Then it follows from Fact \ref{taillaw} that for $n$ large enough and for any particle $u\in\mathcal{Z}[A]$,
\begin{equation}
\Phi_u(z,n)\geq c_{13}/2 (V(u)-z)e^{z-V(u)}\geq \frac{c_{13}}{4} V(u)e^{z-V(u)}.
\end{equation}
(\ref{latticetight}) hence becomes that
\begin{eqnarray*}
\limsup_{n\rightarrow\infty}\P(I_n\geq a_n(0)+ z)&\leq& \E\Big[\prod_{u\in\mathcal{Z}[A]}(1-\frac{c_{13}}{4} V(u)e^{z-V(u)}); \mathcal{Y}_A\Big]+\varepsilon\\
&\leq&\E\Big[\exp\Big(-\frac{c_{13}}{4}e^{z}\sum_{u\in\mathcal{Z}[A]}V(u)e^{-V(u)}\Big)\Big]+\varepsilon.
\end{eqnarray*}
It has been proved that as $A$ goes to infinity, $\sum_{u\in\mathcal{Z}[A]}V(u)e^{-V(u)}$ converges almost surely to some limit $D_\infty$, which is strictly positive on the set of non-extinction of $\mathbb{T}$, (see (5.2) in \cite{elie}).
We end up with
\begin{equation}
\limsup_{z\rightarrow\infty}\limsup_{n\rightarrow\infty}\P(I_n\geq a_n(0)+ z)\leq\varepsilon,
\end{equation}
which completes the proof of Theorem \ref{mainthm}. $\square$

\textbf{Acknowledgments}

I am grateful to my Ph.D. advisor Zhan Shi for his advice and encouragement.

\bigskip
\bigskip



\end{document}